\documentclass{article}
\usepackage{amsmath}
\usepackage{amssymb}
\usepackage{amssymb}
\usepackage{amsmath}
\usepackage{epsfig}
\usepackage{graphicx}
\usepackage{color}
\usepackage[utf8]{inputenc}
\usepackage[english]{babel}
\usepackage{caption}
\usepackage{graphicx}
\usepackage{float}
\usepackage{amsfonts}

\input{tcilatex}
\begin{document}

\begin{center}
\textbf{A Robust Septic Hermite Collocation technique for Dirichlet Boundary
Condition Heat Conduction Equation }

\textbf{\textrm{Sel\c{c}uk KUTLUAY}}$^{\mathbf{\mathrm{a}}}$\textbf{\textrm{%
, Nuri Murat YA\u{G}MURLU}}*$^{\mathbf{\mathrm{a}}}${, }\textbf{\textrm{Ali
Sercan KARAKA\c{S}}}$^{\mathrm{a}}$\textbf{\textrm{{\textrm{\\[0pt]
}}}}$^{\mathbf{\mathrm{{\mathrm{a}}}}}$\textbf{\textrm{{\textbf{\textrm{%
\textbf{In\"{o}n\"{u} University, }}Department of Mathematics, Malatya,
44280, TURKEY.\\[0pt]
}}}}

\textbf{\textrm{{e-mail: \textbf{\textrm{selcuk.kutluay@\textbf{\textrm{%
inonu.edu.tr}}}} }}}

\textbf{\textrm{{e-mail: \textbf{\textrm{murat.yagmurlu@\textbf{\textrm{%
inonu.edu.tr}}}} }}}

\textbf{\textrm{e-mail: 36203614003@ogr.inonu.edu.tr}}
\end{center}

\section{Abstract}

In the present manuscript, approximate solution for 1D heat conduction
equation will be sought with the Septic Hermite Collocation Method (SHCM).
To achieve this goal, by means of the roots of both Chebyschev and Legendre
polinomials used at the inner collocation points, the pseudo code of this
method is found out and applied using Matlab, one of the widely used
symbolic programming platforms. Furthermore, to illustrate the accuracy and
effectiveness of this newly presented scheme, a comparison among analytical
and numerical values is investigated. It has been illustrated that this
scheme is both accurate and effective one and at the same time can be
utilized in a successful way for finding out numerical solutions of several
problems both linear and nonlinear.

\textbf{Keywords:} 1D Heat Equation, Approximate Solutions, Septic Hermite
Collocation Method (SHCM), Finite Elements Method (FEM).

\textbf{AMS classification:} 65M60, 65M12, 65M15

\section{Introduction}

This manuscript is going to dwell on the following 1D heat equation
\begin{equation}
u_{t}=\alpha ^{2}u_{xx},\text{ \ \ }t>0\text{...}x\in \text{ }%
[x_{l},x_{r}],\ \   \label{1}
\end{equation}%
with the following condition given at the beginning of time
\begin{equation}
u(x,t_{initial})=f(x)  \label{2}
\end{equation}%
and conditions given at the boundaries of the solution domain%
\begin{equation}
u(x_{l},t)=u(x_{r},t)=0,\ \ \ \text{ }t>0  \label{3}
\end{equation}%
where $t_{initial}=0,$ $x_{l}=0$, $x_{r}=L,$ $\alpha $ stands for the rod
thermal diffusivity and the $f(x)$ stands for a previously defined function.
Many studies in the literature \cite{1,2,3,4} clearly states that this
initial and boundary value problem (IBVP) can be considered as one of the
most outstanding Partial Differential Equations (PDEs). The equation can
mostly be encountered in various fields of applied sciences and mathematics
for describing the temperature change, or in other words, heat the
distribution, throughout a given solution domain for a given temporal domain.

The manuscript is going to take into consideration the flow of the heat in
1D insulated throughly except with the both of ends of the investigated rod.
The found solutions of this manuscript are going to be presented in terms of
functions given along the rod in spatial and temporal direction. One can
obviously see that this equation has got a vital importance in various areas
of scientific world such as chemistry, physics and mechanics etc. To tell
the truth, the heat equation can be seen as a prototype for parabolic PDE in
several fields of science. Due to its widely utilization and vital place,
both the numerical and exact solutions of those types of equations have
taken attention and essentially have become important in order to
investigate various natural phenomena. Eq. (\ref{1}) is at the same time has
been handled for economical problems in order to simulate several models
with different choices \cite{hikmetmehmetnazan}. The heat equation can at
the same time be described as the flow of heat given in a rod having
diffusion $\alpha ^{2}u_{xx}$ through the rod through having the parameter $%
\alpha $ representing the thermal diffusivity in this rod and $L$ stands for
the length of the considered rod \cite{5}. Kaskar has stated that
electrodynamics, fluid flow, elasticity, electrostatics and other natural
phenomena rely mostly on these equations \cite{kaskar}. One can see that
this problem is considered both analytically and numerically for many years
by many reseachers. Despite of this, the heat equation is still an
outstanding one because of the fact that several natural phenomena can be
stated as PDEs given together with appropriate conditions given at the
initial time and the boundaries of the solution domain. \c{C}a\u{g}lar
\textit{et al.} \cite{hikmetmehmetnazan} have focused on a novel method of
solution of the 1D heat equation by utilizing the third degree B-spline
functions. Dhawan and Kumar \cite{4} have studied temperature change
utilizing the cubic B-spline FEM and gave satisfactory results.
Suarez-Carreno and Rosales-Romero \cite{carrenoromero} have tried to check
the previously obtained results for the heat equation, utilizing both
explicit and implicit methods which are stable for time steps much larger
than those normally utilized in their manuscript by time dependency. The
temporal steps that they have achieved in their work are, for similiar
accuracy, of the same order as those utilized in the implicit methods. Goh
\textit{et al}. \cite{gohmajidismail} have considered a combination of cubic
B-spline method and the finite difference.approach in order to obtain the
solution of 1D advection-diffusion and heat conduction equations. Forward
finite difference approach has been utilized to discretize the temporal
derivative, while cubic B-spline has been applied in order to interpolate
the solution at the time $t$. Lozada-Cruz \textit{et al}.\cite%
{cruzmercedesribeiro} have proved several results on the well-posedness of
the 1D heat equation ans associated stationary problem. and then presented
the application of the Crank-Nicolson for the numer\c{s}ical solution of
considered problem involving the Robin boundary condition. Hooshmandasl
\textit{et al}. \cite{hooshmandasl} have developed Chebyshev wavelet method
together with operational matrices of integration for solution the 1D heat
conduction equation given together with Dirichlet boundary conditions being
fast, simple in mathematical manner and at the same time guarantees the
necessary accuracy for a relative small number of grid points. Han and Dai
\cite{hermiteek4} \ considered two and higher compact finite difference
schemes to obtain solutions for heat conduction equations in 1D with two
test problems: the first one is given together with the Dirichlet boundary
condition while the second one is given together with the Neumann boundary
condition. Kutluay \textit{et al}. \cite{newtrends} \ proposed an effective
numerical approach using cubic Hermite B-spline Collocation method in order
to solve the 1D heat equation. Sun and Zhang \cite{sonek1} proposed a class
of new finite difference scheme, CBVM, in order to solve the 1D heat
equations. They also have shown that CBVM schemes are of high-order accuracy
and unconditionally stable. Patel and Pandya \cite{sonek2} have solved 1D
heat equation given by appropriate Neumann and Dirichlet types
initial-boundary conditions and have utilized a spline collocation method.
Tarmizi \textit{et al}. \cite{sonek3} have used the Spectral method and the
Crank-Nicolson method in order to solve one-dimensional heat equation
together with Dirichlet boundary conditions. Yosaf \textit{et al.} \cite%
{hermiteek3} have provided the implicit compact finite difference method
since it gives a way which is more accurate in order to approximate the
spatial derivative when compared with the explicit finite difference method.

This article will try to find the numerical solution of the 1D heat
conduction problem presented in Eqs.(\ref{1})-(\ref{3}). To accomplish the
goal, Septic Hermite basis functions are going to be utilized. If one
invesitages its historical development, it can be seen that in 1946, first
of all, Schoenberg initiated the basic theory of classical B-splines \cite{4}%
. If one uses the iterative process, linear, quadratic, cubic, quartic etc.
B-spline functions can be easily found. Cox and de Boor have presented the
iterative relationship for calculating the parameters of the new B-spline
functions \cite{8,9} and therefore this relationship is widely known by
their names.Recent years have witnessed a boom both in the theory of spline
function and their applications for obtaining numerical solutions of the the
differential equations approx imately such as those in Refs. \cite%
{hikmetmehmetnazan,14}. Furthermore, various interesting problems are
studied using different FEM methods such as Petrov-Galerkin, Galerkin, least
square, subdomain and collocation method by different odd-numbered degrees
B-spline functions \cite{15,16,17}.

Many differential equations have been solved approximately by various
methods, schemes and techniques using both the usual, trigonometric and
Hermite Septic spline basis functions. Although these methods have some
drawbacks, they are worth to utilize in the application of numerical
methods. One can also say their adavantages over other methods such as (1)
one obtains a diagonal matrix and its storage and manupilation in digital
computers are easy; (2) the computer speed and digital storage area are in
favor of programming.

This article will obtain approximate solutions of the 1D heat conduction
equation given by (\ref{1}) together with the initial (\ref{2}) and boundary
conditions (\ref{3}) using the SHCM, and then those solutions are going to
be compared with the analytical ones \cite{gohmajidismail}. In order to
construct the SHCM scheme, the finite element collocation method\ will be
utilized as in Refs. \cite{20,21}. The second section puts forward the
problem, the method and at the same time presents some useful information
related to the application of the collocation FEM using Septic Hermite
splines. The Section three presents the SHCM for spatial discretization. The
section four proposes the application of the the method in the temporal
direction. The fifth section presents the approximate solutions in tabular
form and at the same time compares the newly obtained solutions with a few
of the previously published ones. The sixth section put forwards a brief
conclusion about the present work.

\section{\textbf{Spatial discretization}}

Throughout this manuscript, the 1D heat conduction equation which has been
widely given by Eq. \ref{1} by the appropriate initial (\ref{2}) and
boundary conditions (\ref{3}) is going to be handled. For this purpose, the
finite element collocation method is selected.as a tool. While applying this
method, the range $[a,b]$ of the problem is discretized in $N$ equal width
finite elements by means of the mesh points $x_{j}$,~~$j=0(1)N$ in such a
way that $x_{l}=x_{0}<x_{1}\cdots <x_{N}=x_{r}$ and $h=x_{j+1}-x_{j}$. To be
more precise, a non-uniform partition of the solution domain might also be
preferred. However, due to the fact that the non-uniform selection would
have increased the digital computer memory storage requirement and at the
same time simulation duration, the uniform partition is usually preferred.
In place of analytical solution $u(x,t)$ a numerical approximation $%
u_{N}(x,t)$ will be used by means of the Septic Hermite basis functions:

\begin{equation}
u_{N}(x,t)=\overset{N+4}{\underset{j=1}{\sum }}a_{j+6k-6}\left( t\right)
H_{ji}\approx u(x,t)  \label{4}
\end{equation}%
in which $a$'s stand for time dependent coefficients , $k$ stands for the
number of elements and finally $i$ stands for the inner-collocation points
used as $i=1,2,3,4,5,6$.

When the following the roots of a six-degree shifted Legendre polinomial
\cite{septicek1}

\begin{eqnarray*}
\xi _{1} &=&0.033765,\quad \xi _{2}=0.169395 \\
\xi _{3} &=&0.380690,\quad \xi _{4}=0.619310 \\
\xi _{4} &=&0.830605,\quad \xi _{6}=0.966235
\end{eqnarray*}

When the following local coordinate system has been utilized over the $%
k^{th} $\ element $\xi =(x-x_{k})/h$ the range $\left[ x_{k},x_{k+1}\right] $
can be changed into $\left[ 0,1\right] $. When this change is used, the
numerical solutions, the first and the second order derivatives of these
solutions, respectively

\begin{align*}
H_{1}\left( \xi \right) & =1-35\xi ^{4}+84\xi ^{5}-70\xi ^{6}+20\xi
^{7},\quad \ \ \ \ \ \ \ \ \ \ \ \ \ H_{2}\left( \xi \right) =h(\xi -20\xi
^{4}+45\xi ^{5}-36\xi ^{6}+10\xi ^{7}) \\
H_{3}\left( \xi \right) & =h^{2}\left( 0.5\xi ^{2}-5\xi ^{4}+10\xi ^{5}-%
\frac{15}{2}\xi ^{6}+2\xi ^{7}\right) ,\quad \text{\ }H_{4}\left( \xi
\right) =h^{3}\left( \frac{1}{6}\xi ^{3}-\frac{2}{3}\xi ^{4}+\xi ^{5}-\frac{2%
}{3}\xi ^{6}+\frac{1}{6}\xi ^{7}\right) \\
H_{5}\left( \xi \right) & =h^{3}\left( -\frac{1}{6}\xi ^{4}+\frac{1}{2}\xi
^{5}-\frac{1}{2}\xi ^{6}+\frac{1}{6}\xi ^{7}\right) ,\quad \ \ \ \ \ \ \ \ \
\ H_{6}\left( \xi \right) =h^{2}\left( \frac{5}{2}\xi ^{4}-7\xi ^{5}+\frac{13%
}{2}\xi ^{6}-2\xi ^{7}\right) \\
H_{7}\left( \xi \right) & =35\xi ^{4}-84\xi ^{5}+70\xi ^{6}-20\xi ^{7},\quad
\ \ \ \ \ \ \ \ \ \ \ \ \ \ \ \ \ \ H_{8}\left( \xi \right) =h(-15\xi
^{4}+39\xi ^{5}-34\xi ^{6}+10\xi ^{7}) \\
& \\
A_{1}\left( \xi \right) & =-140\xi ^{3}+420\xi ^{4}-420\xi ^{5}+140\xi
^{6},\quad \text{\ \ \ \ \ \ \ \ \ }A_{2}\left( \xi \right) =h\left( 1-80\xi
^{3}+225\xi ^{4}-216\xi ^{5}+70\xi ^{6}\right) \\
A_{3}\left( \xi \right) & =h^{2}(\xi -20\xi ^{3}+50\xi ^{4}-45\xi ^{5}+14\xi
^{6}),\quad \ \ \ \ \ \ \ A_{4}\left( \xi \right) =h^{3}\left( \frac{1}{2}%
\xi ^{2}-\frac{8}{3}\xi ^{3}+5\xi ^{4}-4\xi ^{5}+\frac{7}{6}\xi ^{6}\right)
\\
A_{5}\left( \xi \right) & =h^{3}\left( -\frac{2}{3}\xi ^{3}+\frac{5}{2}\xi
^{4}-3\xi ^{5}+\frac{7}{6}\xi ^{6}\right) ,\quad \text{\ \ \ \ \ \ \ \ \ \ \
}A_{6}\left( \xi \right) =h^{2}(10\xi ^{3}-35\xi ^{4}+39\xi ^{5}-14\xi ^{6})
\\
A_{7}\left( \xi \right) & =140\xi ^{3}-420\xi ^{4}+420\xi ^{5}-140\xi
^{6},\quad \text{\ \ \ \ \ \ \ \ \ \ \ }A_{8}\left( \xi \right) =h\left(
-60\xi ^{3}+195\xi ^{4}-204\xi ^{5}+70\xi ^{6}\right) \\
& \\
B_{1}\left( \xi \right) & =-420\xi ^{2}+1680\xi ^{3}-2100\xi ^{4}+840\xi
^{5},\ \ \ \ \ \ \ \ \ B_{2}\left( \xi \right) =h\left( -240\xi ^{2}+900\xi
^{3}+1080\xi ^{4}+420\xi ^{5}\right) \\
B_{3}\left( \xi \right) & =h^{2}\left( 1-60\xi ^{2}+200\xi ^{3}-225\xi
^{4}+84\xi ^{5}\right) ,\ \ \ \ \ B_{4}\left( \xi \right) =h^{3}\left( \xi
-8\xi ^{2}+20\xi ^{3}-20\xi ^{4}+7\xi ^{5}\right) \\
B_{5}\left( \xi \right) & =h^{3}\left( -2\xi ^{2}+10\xi ^{3}-15\xi ^{4}+7\xi
^{5}\right) ,\quad \ \ \ \ \ \ \ \ \ \ \ B_{6}\left( \xi \right)
=h^{2}\left( 30\xi ^{2}-140\xi ^{3}+195\xi ^{4}-84\xi ^{5}\right) \\
B_{7}\left( \xi \right) & =420\xi ^{2}-1680\xi ^{3}+2100\xi ^{4}-840\xi
^{5},\ \ \ \ \ \ \ \ \ \ \ B_{8}\left( \xi \right) =h\left( -180\xi
^{2}+780\xi ^{3}-1020\xi ^{4}+420\xi ^{5}\right) \\
&
\end{align*}%
are obtained. Here $A_{i}$ and $B_{i}$ for $i=1(1)6$ are the first and the
second order derivatives of Septic Hermite functions. Thus, the numerical
approximations are given as follows

\begin{equation}
\left.
\begin{array}{c}
\begin{array}{c}
\text{ \ \ \ \ }%
u_{i}=a_{6k-5}H_{1i}+a_{6k-4}H_{2i}+a_{6k-3}H_{3i}+a_{6k-2}H_{4i}+ \\
a_{6k-1}H_{5i}+a_{6k}H_{6i}+a_{6k+1}H_{7i}+a_{6k+2}H_{8i} \\
\text{ \ }hu_{i}^{\prime
}=a_{6k-5}A_{1i}+a_{6k-4}A_{2i}+a_{6k-3}A_{3i}+a_{6k-2}A_{4i}+ \\
a_{6k-1}A_{5i}+a_{6k}A_{6i}+a_{6k+1}A_{7i}+a_{6k+2}A_{8i}%
\end{array}
\\
h^{2}u_{i}=a_{6k-5}B_{1i}+a_{6k-4}B_{2i}+a_{6k-3}B_{3i}+a_{6k-2}B_{4i}+ \\
a_{6k-1}B_{5i}+a_{6k}B_{6i}+a_{6k+1}B_{7i}+a_{6k+2}B_{8i}%
\end{array}%
\right\} ,\quad i=1,2,3,4,5,6.  \label{5}
\end{equation}

While applying the process, on the spatial direction Hermite Septic spline
basis functions and on the temporal direction the forward finite difference
formula will be used. This choice has the advantages of several vital
characteristics such as easy to handle algorithms they produce and low level
of storage requirement. Besides those advantages, both of the non-linear and
linear systems resulted from the usage of splines are usually not
ill-conditioned and thus permit the required coefficients to be found out in
a easy manner. Furthermore, the newly obtained approximate solutions
generally don't result in numerical instability. The last but not the least,
the matrix systems resulting from the usage of splines are usually sparsly
loaded band matrixes and allow easy implementation in digital computers \cite%
{15}.

\section{\textbf{Temporal discretization}}

When the aforementioned 1D heat conduction equation given by Eq. (\ref{1})
as follows is discretized
\begin{equation*}
u_{t}-\alpha ^{2}u_{xx}=0.
\end{equation*}%
one can use Crank-Nicolson type formula. Firstly, one discretizes Eq. (\ref%
{1}) as%
\begin{equation*}
\frac{u^{n+1}-u^{n}}{\Delta t}-\alpha ^{2}\left[ \frac{%
(u_{xx})^{n+1}+(u_{xx})^{n}}{2}\right] =0.
\end{equation*}%
If one seperates the variables such that the next time level ones are placed
on the left hand side and the ones for previous time level are placed on the
right hand side, then
\begin{equation}
\frac{u^{n+1}}{\Delta t}-\alpha ^{2}\frac{(u_{xx})^{n+1}}{2}=\frac{u^{n}}{%
\Delta t}+\alpha ^{2}\frac{(u_{xx})^{n}}{2}  \label{6}
\end{equation}%
is obtained.

If one substituties Eq. (\ref{5}) in Eq. (\ref{6}), the the following
difference equation system is found out for the coefficients $\mathbf{%
\mathbf{a}}$ with $4N$ difference equations and $4N+2$ coeffiicents

\begin{equation}
\begin{array}{l}
\text{ \ }\frac{1}{\Delta t}\left[
\begin{array}{c}
a_{6k-5}^{n+1}H_{1i}+a_{6k-4}^{n+1}H_{2i}+a_{6k-3}^{n+1}H_{3i}+a_{6k-2}^{n+1}H_{4i}+
\\
a_{6k-1}^{n+1}H_{5i}+a_{6k}^{n+1}H_{6i}+a_{6k+1}^{n+1}H_{7i}+a_{6k+2}^{n+1}H_{8i}%
\end{array}%
\right] \\
\\
-\frac{\alpha ^{2}}{2h^{2}}\left[
\begin{array}{c}
a_{6k-5}^{n+1}B_{1i}+a_{6k-4}^{n+1}B_{2i}+a_{6k-3}^{n+1}B_{3i}+a_{6k-2}^{n+1}B_{4i}+
\\
a_{6k-1}^{n+1}B_{5i}+a_{6k}^{n+1}B_{6i}+a_{6k+1}^{n+1}B_{7i}+a_{6k+2}^{n+1}B_{8i}%
\end{array}%
\right] \\
\\
=\frac{1}{\Delta t}\left[
\begin{array}{c}
a_{6k-5}^{n}H_{1i}+a_{6k-4}^{n}H_{2i}+a_{6k-3}^{n}H_{3i}+a_{6k-2}^{n}H_{4i}+
\\
a_{6k-1}^{n}H_{5i}+a_{6k}^{n}H_{6i}+a_{6k+1}^{n}H_{7i}+a_{6k+2}^{n}H_{8i}%
\end{array}%
\right] \\
\\
+\frac{\alpha ^{2}}{2h^{2}}\left[
\begin{array}{c}
a_{6k-5}^{n}B_{1i}+a_{6k-4}^{n}B_{2i}+a_{6k-3}^{n}B_{3i}+a_{6k-2}^{n}B_{4i}+
\\
a_{6k-1}^{n}B_{5i}+a_{6k}^{n}B_{6i}+a_{6k+1}^{n}B_{7i}+a_{6k+2}^{n}B_{8i}%
\end{array}%
\right]%
\end{array}
\label{7}
\end{equation}

Those newly obtained equations are recursive in nature. Thus the element
coefficients vector $\mathbf{a}%
^{n}=(a_{1}^{n},...,a_{6N+1}^{n},a_{6N+2}^{n}) $ can be found in a recursive
manner, where $t_{n}=n(\Delta t),$ $n=1(1)M$ till the desired time $T$. If
one utilizes the conditions presented at the boundary of the solution domain
by Eq.(\ref{3}) and eliminates the coefficients $a_{1}^{n},a_{6N+1}^{n}$ in
Eq. (\ref{7}) as follows: Using the boundary condition given at the left of
the solution domain $%
u(x_{l},t)=a_{1}^{n}H_{11}+a_{2}^{n}H_{21}+a_{3}^{n}H_{31}+a_{4}^{n}H_{41}+a_{5}^{n}H_{51}+a_{6}^{n}H_{61}+a_{7}^{n}H_{71}+a_{8}^{n}H_{81}=0,
$ since $H_{m1}=0$ for $m=2(1)8$ and $H_{11}\neq 0,$ the condition $%
a_{1}^{n}=0$ is obtained. In a similar way, using the right boundary
condition $%
u(x_{r},t)=a_{6N-5}^{n}H_{16}+a_{6N-4}^{n}H_{26}+a_{6N-3}^{n}H_{36}+a_{6N-2}^{n}H_{46}+a_{6N-1}^{n}H_{56}+a_{6N}^{n}H_{66}+a_{6N+1}^{n}H_{76}+a_{6N+2}^{n}H_{86}=0,
$ since $H_{m6}=0$ for $m=1,2,3,4,5,6,8$ and $H_{76}\neq 0$, the condition $%
a_{6N+1}^{n}=0$ is found. Finally, a system which is solvable is obtained as
follows

\begin{equation}
\mathbf{La}^{n+1}=\mathbf{Ra}^{n}.  \label{recursive}
\end{equation}%
In this iterative formulae, the matrices $\mathbf{\mathbf{L}}$ and $\mathbf{R%
}$ are square $6N\times 6N$ diagonal band matrices. The matrices $\mathbf{a}%
^{n+1}$ and $\mathbf{a}^{n}$ are $6N$ column matrices.

The coefficients $\mathbf{a}_{i}$ $(i=1(1)6N)$ found using the system of
equations in Eq.(\ref{recursive}) have been obtained and the numerical
solutions of the present equation at the next temporal levels are
calculated. The iterative procedure is made in a successive manner for $%
t_{n}=n\Delta t$ $(n=1(1)M)$ until the desired time $T$. To be able to start
this iterative calculation, the initial\ vector\ $\mathbf{a}^{0}$ having the
entries $\mathbf{a}_{i0}$ $(i=1(1)6N)$ should be found. This initial vector
is found out with the usage of the initial condition given together with the
governing equation. When one substitutites these values in their places, Eq.
(\ref{7}) results in $6N$ unknowns. To obtain the solution of this system,
an algorithm written in MatlabR2021a is run on a digital computer.

\subsection{\textbf{The initial state}}

To be able to initiate iteration process, one needs the vector at the
initial time. The vector $\mathbf{a}^{0}$ is found out using the conditions
given at the initial time and at the boundaries of the solution domain.
Thus, the numerical approximation given by Eq. ($\ref{4})$ should be
rewritten in an appropriate way for the initial condition as follows

\begin{equation*}
u_{N}(x,t)=\overset{N+4}{\underset{j=1}{\sum }}a_{j+6k-6}^{0}\left( t\right)
H_{ji}\approx u(x,t)
\end{equation*}%
in which the $a_{m}^{0}$'s are unknown coefficients to be calculated. One
requires that the initial approximate solution $u_{N}(x,0)$ satisfies the
initial analytical solution $u(x,0)$ as follows

\begin{equation*}
u_{N}(x_{i},0)=u(x_{i},0),\ \ \ \ \ \ i=0(1)N
\end{equation*}%
Finally, those calculations result in the matrix equation of the following
form {\small
\begin{equation}
\mathbf{Wa}^{0}=\mathbf{b}  \label{denkek}
\end{equation}%
} where {\small {\
\begin{equation*}
\mathbf{W}=\left[
\begin{array}{ccccccccccccccc}
H_{21} & H_{31} & H_{41} & H_{51} & H_{61} & H_{71} & H_{81} &  &  &  &  &
&  &  &  \\
H_{22} & H_{32} & H_{42} & H_{52} & H_{62} & H_{72} & H_{82} &  &  &  &  &
&  &  &  \\
H_{23} & H_{33} & H_{43} & H_{53} & H_{63} & H_{73} & H_{83} &  &  &  &  &
&  &  &  \\
H_{24} & H_{34} & H_{44} & H_{54} & H_{64} & H_{74} & H_{84} &  &  &  &  &
&  &  &  \\
H_{25} & H_{35} & H_{45} & H_{55} & H_{65} & H_{75} & H_{85} &  &  &  &  &
&  &  &  \\
H_{26} & H_{36} & H_{46} & H_{56} & H_{66} & H_{76} & H_{86} &  &  &  &  &
&  &  &  \\
&  &  &  &  & H_{11} & H_{21} & H_{31} & H_{41} & H_{51} & H_{61} & H_{71} &
H_{81} &  &  \\
&  &  &  &  & H_{12} & H_{22} & H_{32} & H_{42} & H_{52} & H_{62} & H_{72} &
H_{82} &  &  \\
&  &  &  &  & H_{13} & H_{23} & H_{33} & H_{43} & H_{53} & H_{63} & H_{73} &
H_{83} &  &  \\
&  &  &  &  & H_{14} & H_{24} & H_{34} & H_{44} & H_{54} & H_{64} & H_{74} &
H_{84} &  &  \\
&  &  &  &  & H_{15} & H_{25} & H_{35} & H_{45} & H_{55} & H_{65} & H_{75} &
H_{85} &  &  \\
&  &  &  &  & H_{16} & H_{26} & H_{36} & H_{46} & H_{56} & H_{66} & H_{76} &
H_{86} &  &  \\
&  &  &  &  &  & \ddots & \ddots & \ddots & \ddots & \ddots & \ddots & \ddots
& \ddots &  \\
&  &  &  &  &  &  &  & H_{11} & H_{21} & H_{31} & H_{41} & H_{51} & H_{61} &
H_{81} \\
&  &  &  &  &  &  &  & H_{12} & H_{22} & H_{32} & H_{42} & H_{52} & H_{62} &
H_{82} \\
&  &  &  &  &  &  &  & H_{13} & H_{23} & H_{33} & H_{43} & H_{53} & H_{63} &
H_{83} \\
&  &  &  &  &  &  &  & H_{14} & H_{24} & H_{34} & H_{44} & H_{54} & H_{64} &
H_{84} \\
&  &  &  &  &  &  &  & H_{15} & H_{25} & H_{35} & H_{45} & H_{55} & H_{65} &
H_{85} \\
&  &  &  &  &  &  &  & H_{16} & H_{26} & H_{36} & H_{46} & H_{56} & H_{66} &
H_{86}%
\end{array}%
\right] ,
\end{equation*}%
}%
\begin{equation*}
\begin{array}{l}
\mathbf{a}{^{0}}=(a_{2},a_{3},a_{4},\ldots ,a_{6N-1},a_{6N},a_{6N+2})^{T}%
\end{array}%
\end{equation*}%
}and {\small \
\begin{equation*}
\begin{array}{l}
\mathbf{b}%
=(u(x_{11},0),u(x_{12},0),u(x_{13},0),u(x_{14},0),u(x_{15},0),u(x_{16},0),%
\ldots ,u(x_{N5},0),u(x_{N6},0))^{T}.%
\end{array}%
\end{equation*}%
}

Therefore, the required initial coeffients to initiate the scheme given by
Eq. (\ref{recursive}) are calculated by Eq. (\ref{denkek}) and next the
iterative procedure is made in a repetitive manner till the final time $T$.

\section{Numerical results}

The present section will present the application of the collocation FEM
using Quintic hermite basis functions to a control problem. For the control
problem, the condition function given at the initial time is going to be
used as $f(x)=\sin (\pi x)$. The closed interval $[0,1]$ will be taken as
the solution domain of the control problem. The analytical solution for the
governing problem is \cite{hermiteek3,hermiteek4}

\begin{equation*}
u(x,t)=\frac{\sin (\pi x)}{e^{\alpha ^{2}\pi ^{2}t}}
\end{equation*}

Due to the fact that the control problem has an analytical solution, the
error norms $L_{2}$ and $L_{\infty }$ formulated as, respectively, will be
used to test the accuracy and validity of this method
\begin{equation*}
L_{2}=\sqrt{\left( h\sum_{i=1}^{N}\left\vert u_{i}-(u_{N})_{i}\right\vert
^{2}\right) },\text{\qquad }L_{\infty }=\max_{1\leq i\leq N}\left\vert
u_{i}-(u_{N})_{i}\right\vert .
\end{equation*}%
This manuscript carries out all of the numerical computations using both
Septic Hermite Collocation Method with Legendre roots (SHCM-L) and Septic
Hermite Collocation Method with Chebyshev roots (SHCM-C). Those calculations
have been made with MATLAB R2021a on Intel (R) Core(TM) i57 3330SU CPU
@2.70Ghz computer having 6 GB of RAM.

Various space step sizes, time step sizes and desired values of final time
are used to control the effectiveness and accuracy of the present method.

\begin{table}[tbp]
\caption{A comparison of $\ $the computed error norm\ $L_{2}$ by the
proposed scheme with those given in Ref.\protect\cite{hermiteek3},
\protect\cite{newtrends} for $\Delta t=1/1000$ and $k=\Delta
t=0.01,0.005,0.0025$ ($t_{final}=1,\protect\alpha =1,$ $0\leq x\leq 1$).}
\label{hermitetbl1}${\scriptsize
\begin{tabular}{ccccccccc}
\hline
&  & \multicolumn{7}{c}{$L_{2}$} \\ \hline
$\Delta t$ &  & SHCM-L &  & SHCM-C &  & \cite{hermiteek3} &  & \cite%
{newtrends}CHCM-L \\ \hline
$0.01$ &  & $7.1591\times 10^{-7}$ &  & $7.1591\times 10^{-7}$ &  & $%
4.2273\times 10^{-4}$ &  & $4.1333\times 10^{-7}$ \\
$0.005$ &  & $1.7931\times 10^{-7}$ &  & $1.7932\times 10^{-7}$ &  & $%
1.0560\times 10^{-4}$ &  & $1.0353\times 10^{-7}$ \\
$0.0025$ &  & $4.4851\times 10^{-8}$ &  & $4.4851\times 10^{-8}$ &  & $%
2.6395\times 10^{-5}$ &  & $2.5895\times 10^{-8}$ \\ \hline
\end{tabular}%
\ \ \ \ \ \ \ \ \ \ \ \ \ }$%
\end{table}

Table \ref{hermitetbl1}, illustrates a clear comparison of the newly
computed error norm values\ $L_{2}$ of this scheme with those given in Ref.%
\cite{hermiteek3}, \cite{newtrends} for values of $k=\Delta
t=1/100,5/1000,25/10000$ and $h=\Delta x=1/1000$ ($0\leq x\leq 1,$ $%
t_{final}=1,\alpha =1$). It is obvios from this table that when those values
of $\Delta t$ decrease, at the same time the error norm\ $L_{2}$ values
decrease. Furthermore, it is seen that the obtained values using both SHCM-C
and SHCM-L clearly more the better than those found by Ref. \cite{hermiteek3}

\begin{table}[tbp]
\caption{A comparison of $\ $the computed error norm\ $L_{2}$ of the
proposed scheme with those given in Ref.\protect\cite{hermiteek4},
\protect\cite{newtrends} for $N=5,10,20$ and $k=\Delta t=1/10^{6}$ ($%
t_{final}=1,\protect\alpha =1,$ $0\leq x\leq 1$).}
\label{hermitetbl2}${\scriptsize
\begin{tabular}{cccccccccc}
\hline
& \multicolumn{9}{c}{$L_{2}$} \\ \hline
$\Delta x$ & SHCM-L &  & SHCM-C &  & \cite{hermiteek4}(CN-I) &  & \cite%
{hermiteek4}(CN-II) &  & \cite{newtrends}CHCM-L \\ \hline
$1/5$ & $1.2153\times 10^{-14}$ &  & $1.0075\times 10^{-12}$ &  &
\multicolumn{1}{l}{$4.7696\times 10^{-6}$} & \multicolumn{1}{l}{} &
\multicolumn{1}{l}{$8.5859\times 10^{-3}$} & \multicolumn{1}{l}{} &
\multicolumn{1}{l}{$3.5716\times 10^{-8}$} \\
$1/10$ & $1.9449\times 10^{-15}$ &  & $2.3323\times 10^{-14}$ &  &
\multicolumn{1}{l}{$7.7143\times 10^{-9}$} & \multicolumn{1}{l}{} &
\multicolumn{1}{l}{$2.1412\times 10^{-3}$} & \multicolumn{1}{l}{} &
\multicolumn{1}{l}{$2.2848\times 10^{-9}$} \\
$1/20$ & $4.6215\times 10^{-15}$ &  & $7.4211\times 10^{-15}$ &  &
\multicolumn{1}{l}{$1.8820\times 10^{-11}$} & \multicolumn{1}{l}{} &
\multicolumn{1}{l}{$5.3498\times 10^{-4}$} & \multicolumn{1}{l}{} &
\multicolumn{1}{l}{$1.4361\times 10^{-10}$} \\
$1/40$ & $8.1433\times 10^{-15}$ &  & $9.6794\times 10^{-15}$ &  &  &  &  &
&  \\ \hline
\end{tabular}%
\ \ \ \ \ \ \ \ \ \ \ \ \ }$%
\end{table}

From Table \ref{hermitetbl2}, a clear comparison of the newly computed error
norm values\ $L_{2}$ of this scheme with those found by Ref.\cite{hermiteek4}%
, \cite{newtrends} for several values of $\Delta x=h=1/5,1/10,1/20$ and $%
k=\Delta t=1/10^{6}$ ($0\leq x\leq 1,\alpha =1,$ $t_{final}=1$) can be seen.
It can be easily seen from the table that when the number of elements is
increased, the values of the error norm $L_{2}$ decrease. Similarly, this
table obviously shows that the results found out by utilizing the present
scheme are in general better than those of compared ones among others as in
Ref. \cite{hermiteek4}.

\begin{table}[tbp]
\caption{A comparison of $\ $the computed error norm\ $L_{2}$ of the
proposed scheme with those given in Ref.\protect\cite{hermiteek3},
\protect\cite{newtrends} for $N=10,20,40$ and $k=\Delta t=1/10^{6}$ ($%
t_{final}=0.1,\protect\alpha =1,$ $0\leq x\leq 1$).}
\label{hermitetbl3}${\scriptsize
\begin{tabular}{ccccccccc}
\hline
$\Delta x$ &  & SHCM-L &  & \multicolumn{3}{c}{\cite{hermiteek3}} &  & \cite%
{newtrends}CHCM-L \\ \hline
&  &  &  & $\theta =0.1$ &  & $\theta =0.2$ &  &  \\
$1/10$ &  & $6.2482\times 10^{-13}$ &  & \multicolumn{1}{l}{$1.6534\times
10^{-4}$} & \multicolumn{1}{l}{} & \multicolumn{1}{l}{$2.4968\times 10^{-4}$}
& \multicolumn{1}{l}{} & \multicolumn{1}{l}{$1.6957\times 10^{-6}$} \\
$1/20$ &  & $3.3280\times 10^{-12}$ &  & \multicolumn{1}{l}{$4.3906\times
10^{-7}$} & \multicolumn{1}{l}{} & \multicolumn{1}{l}{$7.8152\times 10^{-7}$}
& \multicolumn{1}{l}{} & \multicolumn{1}{l}{$1.0426\times 10^{-7}$} \\
$1/40$ &  & $5.8595\times 10^{-12}$ &  & \multicolumn{1}{l}{$8.6154\times
10^{-10}$} & \multicolumn{1}{l}{} & \multicolumn{1}{l}{$8.0111\times
10^{-10} $} & \multicolumn{1}{l}{} & \multicolumn{1}{l}{$6.4916\times
10^{-9} $} \\ \hline
&  &  &  & $\theta =0.3$ &  & $\theta =0.4$ &  &  \\
$1/10$ &  & $6.2482\times 10^{-13}$ &  & \multicolumn{1}{l}{$3.2357\times
10^{-4}$} & \multicolumn{1}{l}{} & \multicolumn{1}{l}{$3.6276\times 10^{-4}$}
& \multicolumn{1}{l}{} & \multicolumn{1}{l}{$1.6957\times 10^{-6}$} \\
$1/20$ &  & $3.3280\times 10^{-12}$ &  & \multicolumn{1}{l}{$1.0791\times
10^{-6}$} & \multicolumn{1}{l}{} & \multicolumn{1}{l}{$1.2556\times 10^{-6}$}
& \multicolumn{1}{l}{} & \multicolumn{1}{l}{$1.0426\times 10^{-7}$} \\
$1/40$ &  & $5.8595\times 10^{-12}$ &  & \multicolumn{1}{l}{$9.2171\times
10^{-10}$} & \multicolumn{1}{l}{} & \multicolumn{1}{l}{$1.1166\times 10^{-9}$%
} & \multicolumn{1}{l}{} & \multicolumn{1}{l}{$6.4916\times 10^{-9}$} \\
\hline
&  &  &  & $\theta =0.5$ &  & $\theta =0.6$ &  &  \\
$1/10$ &  & $6.2482\times 10^{-13}$ &  & \multicolumn{1}{l}{$3.5225\times
10^{-4}$} & \multicolumn{1}{l}{} & \multicolumn{1}{l}{$2.8316\times 10^{-4}$}
& \multicolumn{1}{l}{} & \multicolumn{1}{l}{$1.6957\times 10^{-6}$} \\
$1/20$ &  & $3.3280\times 10^{-12}$ &  & \multicolumn{1}{l}{$1.2494\times
10^{-6}$} & \multicolumn{1}{l}{} & \multicolumn{1}{l}{$1.0109\times 10^{-6}$}
& \multicolumn{1}{l}{} & \multicolumn{1}{l}{$1.0426\times 10^{-7}$} \\
$1/40$ &  & $5.8595\times 10^{-12}$ &  & \multicolumn{1}{l}{$1.2167\times
10^{-9}$} & \multicolumn{1}{l}{} & \multicolumn{1}{l}{$1.1107\times 10^{-9}$}
& \multicolumn{1}{l}{} & \multicolumn{1}{l}{$6.4916\times 10^{-9}$} \\ \hline
&  &  &  & $\theta =0.7$ &  & $\theta =0.8$ &  &  \\
$1/10$ &  & $6.2482\times 10^{-13}$ &  & \multicolumn{1}{l}{$1.5954\times
10^{-4}$} & \multicolumn{1}{l}{} & \multicolumn{1}{l}{$1.9640\times 10^{-4}$}
& \multicolumn{1}{l}{} & \multicolumn{1}{l}{$1.6957\times 10^{-6}$} \\
$1/20$ &  & $3.3280\times 10^{-12}$ &  & \multicolumn{1}{l}{$5.4419\times
10^{-7}$} & \multicolumn{1}{l}{} & \multicolumn{1}{l}{$9.1503\times 10^{-7}$}
& \multicolumn{1}{l}{} & \multicolumn{1}{l}{$1.0426\times 10^{-7}$} \\
$1/40$ &  & $5.8595\times 10^{-12}$ &  & \multicolumn{1}{l}{$7.4804\times
10^{-10}$} & \multicolumn{1}{l}{} & \multicolumn{1}{l}{$8.6751\times
10^{-10} $} & \multicolumn{1}{l}{} & \multicolumn{1}{l}{$6.4916\times
10^{-9} $} \\ \hline
&  &  &  & $\theta =0.9$ &  & $\theta =0.95$ &  &  \\
$1/10$ &  & $6.2482\times 10^{-13}$ &  & \multicolumn{1}{l}{$4.5197\times
10^{-4}$} & \multicolumn{1}{l}{} & \multicolumn{1}{l}{$6.1645\times 10^{-4}$}
& \multicolumn{1}{l}{} & \multicolumn{1}{l}{$1.6957\times 10^{-6}$} \\
$1/20$ &  & $3.3280\times 10^{-12}$ &  & \multicolumn{1}{l}{$2.0820\times
10^{-6}$} & \multicolumn{1}{l}{} & \multicolumn{1}{l}{$2.8818\times 10^{-6}$}
& \multicolumn{1}{l}{} & \multicolumn{1}{l}{$1.0426\times 10^{-7}$} \\
$1/40$ &  & $5.8595\times 10^{-12}$ &  & \multicolumn{1}{l}{$2.4668\times
10^{-9}$} & \multicolumn{1}{l}{} & \multicolumn{1}{l}{$3.6858\times 10^{-9}$}
& \multicolumn{1}{l}{} & \multicolumn{1}{l}{$6.4916\times 10^{-9}$} \\
\hline\cline{2-9}
\end{tabular}%
\ \ \ \ \ \ \ \ \ \ \ \ \ }$%
\end{table}

One can clearly see from Table \ref{hermitetbl3} that a clear comparison of
the computed error norm values\ $L_{2}$ of this scheme with these given by
Ref.\cite{hermiteek3}, \cite{newtrends} for $h=\Delta x=1/10,1/20,1/40$ and $%
k=\Delta t=1/10^{6}$ ($0\leq x\leq 1,\alpha =1,$ $t_{final}=0.1$). It can be
obviously seen that the new results found utilizing the present scheme are
much better than the results in Ref. \cite{hermiteek3}.

\begin{table}[tbp]
\caption{A comparison of $\ $the computed error norm\ $L_{\infty }$ of the
proposed scheme with those given in Refs. \protect\cite{gohmajidismail},
\protect\cite{newtrends}and \protect\cite{sonek1} for various values of $%
h=\Delta x=k=\Delta t$ ($t_{final}=1,\protect\alpha =1,$ $0\leq x\leq 1,$ ).}
\label{hermitetbl4}
\begin{center}
${\scriptsize
\begin{tabular}{ccccccc}
\hline
$h=k$ & SHCM-L & SHCM-C & \cite{newtrends}CHCM-L & \cite{gohmajidismail} &
\cite{sonek1}(CN) & \cite{sonek1}(CBVM) \\ \hline
\multicolumn{1}{l}{$0.2$} & $5.1578\times 10^{-5}$ & $5.1552\times 10^{-5}$
& $5.8498\times 10^{-5}$ & \multicolumn{1}{l}{$1.4145\times 10^{-1}$} &
\multicolumn{1}{l}{$1.1\times 10^{-1}$} & \multicolumn{1}{l}{$2.8\times
10^{-2}$} \\
\multicolumn{1}{l}{$0.1$} & $3.1586\times 10^{-5}$ & $3.1587\times 10^{-5}$
& $3.1584\times 10^{-5}$ & \multicolumn{1}{l}{$3.7195\times 10^{-2}$} &
\multicolumn{1}{l}{$3.0\times 10^{-2}$} & \multicolumn{1}{l}{$3.8\times
10^{-3}$} \\
\multicolumn{1}{l}{$0.05$} & $9.7106\times 10^{-6}$ & $9.7107\times 10^{-6}$
& $9.7065\times 10^{-6}$ & \multicolumn{1}{l}{$8.4588\times 10^{-3}$} &
\multicolumn{1}{l}{$6.9\times 10^{-3}$} & \multicolumn{1}{l}{$2.7\times
10^{-4}$} \\
\multicolumn{1}{l}{$0.025$} & $2.5489\times 10^{-6}$ & $2.5489\times 10^{-6}$
& $2.5485\times 10^{-6}$ & \multicolumn{1}{l}{$2.0698\times 10^{-3}$} &
\multicolumn{1}{l}{$1.7\times 10^{-3}$} & \multicolumn{1}{l}{$1.3\times
10^{-5}$} \\
\multicolumn{1}{l}{$0.0125$} & $6.4490\times 10^{-7}$ & $6.4490\times
10^{-7} $ & $6.4488\times 10^{-7}$ & \multicolumn{1}{l}{$5.1473\times
10^{-4} $} & \multicolumn{1}{l}{$4.2\times 10^{-4}$} & \multicolumn{1}{l}{$%
5.1\times 10^{-7}$} \\
\multicolumn{1}{l}{$0.00625$} & $1.6171\times 10^{-7}$ & $1.6171\times
10^{-7}$ & $1.6171\times 10^{-7}$ & \multicolumn{1}{l}{-} &
\multicolumn{1}{l}{$1.1\times 10^{-4}$} & \multicolumn{1}{l}{$3.6\times
10^{-8}$} \\
\multicolumn{1}{l}{$0.01$} & $4.1333\times 10^{-7}$ & $4.1333\times 10^{-7}$
& $4.1332\times 10^{-7}$ &  &  &  \\
\multicolumn{1}{l}{$0.005$} & $1.0353\times 10^{-7}$ & $1.0353\times 10^{-7}$
& $1.0353\times 10^{-7}$ &  &  &  \\
\multicolumn{1}{l}{$0.0025$} & $2.5895\times 10^{-8}$ & $2.5895\times
10^{-8} $ & $2.5895\times 10^{-8}$ &  &  &  \\
\multicolumn{1}{l}{$0.002$} & $1.6574\times 10^{-8}$ & $1.6574\times 10^{-8}$
& $1.6574\times 10^{-8}$ &  &  &  \\
\multicolumn{1}{l}{$0.001$} & $4.1437\times 10^{-9}$ & $4.1437\times 10^{-9}$
& $4.1437\times 10^{-9}$ &  &  &  \\ \hline
\end{tabular}%
\ \ \ \ \ \ \ \ \ \ \ \ }$%
\end{center}
\end{table}

One can see from Table \ref{hermitetbl4} a clear comparison of the error
norm values $L_{\infty }$ of this scheme with those given by Refs. \cite%
{gohmajidismail}, \cite{newtrends} and \cite{sonek1} for various values of $%
h=\Delta x=k=\Delta t$ ($\alpha =1,$ $0\leq x\leq 1,$ $t_{final}=1$). It can
be obviously umdesrtood from this table that the approximations are better
or in good agreement with those of the compared ones.

\begin{table}[tbp]
\caption{A comparison of $\ $the computed error norms\ $L_{2}$ and $%
L_{\infty }$ of the proposed scheme with those given in Ref. \protect\cite%
{hooshmandasl}, \protect\cite{newtrends} for $\Delta =1/16,$ $k=\Delta
t=0.01 $ at various $t_{final}$ values ($0\leq x\leq 1,\protect\alpha =1$). }
\label{hermitetbl5}
\begin{center}
${\scriptsize
\begin{tabular}{cccccccc}
\hline
$t_{final}$ & \multicolumn{3}{c}{$L_{2}$} &  & \multicolumn{3}{c}{$L_{\infty
}$} \\ \cline{1-4}\cline{6-8}
& SHCM-L & \cite{newtrends}CHCM-L & \cite{hooshmandasl} &  & SHCM-L & \cite%
{newtrends}CHCM-L & \cite{hooshmandasl} \\ \hline
$0.1$ & $5.1774\times 10^{-4}$ & $2.9917\times 10^{-4}$ & $4.86\times
10^{-3} $ &  & $2.9891\times 10^{-4}$ & $2.9891\times 10^{-4}$ & $6.79\times
10^{-3}$ \\
$0.3$ & $2.1558\times 10^{-4}$ & $1.2457\times 10^{-4}$ & $8.87\times
10^{-5} $ &  & $1.2447\times 10^{-4}$ & $1.2447\times 10^{-4}$ & $3.76\times
10^{-4}$ \\
$0.5$ & $4.9872\times 10^{-5}$ & $2.8818\times 10^{-5}$ & $1.73\times
10^{-3} $ &  & $2.8793\times 10^{-5}$ & $2.8793\times 10^{-5}$ & $2.44\times
10^{-4}$ \\
$0.7$ & $9.6911\times 10^{-6}$ & $5.5999\times 10^{-6}$ & $2.04\times
10^{-4} $ &  & $5.5950\times 10^{-6}$ & $5.5953\times 10^{-6}$ & $3.17\times
10^{-4}$ \\
$0.9$ & $1.7294\times 10^{-6}$ & $9.9934\times 10^{-7}$ & $2.14\times
10^{-3} $ &  & $9.9847\times 10^{-7}$ & $9.9856\times 10^{-7}$ & $3.14\times
10^{-3}$ \\
$1.0$ & $7.1591\times 10^{-7}$ & $4.1368\times 10^{-7}$ & $2.15\times
10^{-3} $ &  & $4.1332\times 10^{-7}$ & $4.1338\times 10^{-7}$ & $3.32\times
10^{-3}$ \\ \hline
\end{tabular}%
\ \ \ \ \ \ \ \ \ \ \ \ \ \ }$%
\end{center}
\end{table}

Finally, one can celarly see in Table \ref{hermitetbl5} a comparison of the
error norm values of $L_{2}$ and $L_{\infty }$ for $N=16,$ $\alpha =1,k=0.01$
with those of Ref.\cite{hooshmandasl}, \cite{newtrends} at various $%
t_{final} $ values. When the table is investigated, it is obvious that the
newly computed results are better or in good agreement with those of
compared ones in view of the error norms\ $L_{2}$ and $L_{\infty }$.

It can be seen that as the step size in spatial dimension is decreased, the
numerical results also approach to their analytical. Tables \ref{hermitetbl2}
and \ref{hermitetbl3}.clearly show this fact.

\section{Conclusion}

This manuscript has successfully applied the collocation finite element
method resulting an implicit linear algebraic system to find out the
numerical solutions for the 1D heat conduction equation. The control problem
has shown that the approximations have been in good harmony with the
analytical solutions. As a conclusion, the newly presented numerical scheme,
with easy implementation, results in efficient and accurate solutions. For a
prospective future study, this method can be successfully utilized to obtain
accurate numerical solutions of various PDEs that has a vital role in
describing natural phenomena which is both linear and nonlinear and widely
encountered in our daily lives.

\section{Declaration of Ethical Standards}

The authors of the present manuscript clearly declare that all of the
methods and schemes used in the manuscript don't need any ethical committee
and/or legal special requirement or permission.

The authors also declare that this manuscript has been presented in "The 7th
International Conference on Computational Mathematics and Engineering
Sciences / 20-21 May. 2023, Elaz\i \u{g} -- T\"{u}rkiye."

\subsection*{Financial disclosure}

One of the authors, Ali Sercan KARAKA\c{S}, is funded by 2211-A National PhD
Scholarship Program, The Scientific and Technological Research Council of T%
\"{U}RK\.{I}YE.


\begin{thebibliography}{99}
\bibitem{1} N. \"{O}z\i \c{s}\i k, Heat Conduction, CRC Press Taylor and
Francis Group, (2017).

\bibitem{2} D.V. Widder, The Heat Equation, Academic Press, (1976).

\bibitem{3} J.R. Canon, The one-Dimensional Heat Equation, Cambridge
University Pres., (1984).

\bibitem{4} S. Dhawan, S. Kumar, A numerical solution of one dimensional
Heat Equation Using Cubic B-spline Basis Functions, \textbf{International
Journal of Research and Reviews in Applied Sciences} , (2009) 71-77.

\bibitem{5} A.M. Wazwaz, Partial Differential Equations methods and
Applications, Saint Xavier University Pres., (2002).

\bibitem{hikmetmehmetnazan} H. \c{C}a\u{g}lar, M. \"{O}zer, N. \c{C}a\u{g}%
lar, The numerical solution of the one-dimensional heat equation by using
third degree B-spline functions, \textbf{Chaos, Solitons and Fractals} 38
(2008) 1197--1201

\bibitem{carrenoromero} F. Suarez-Carreno and L. Rosales-Romero, Convergency
and Stability of Explicit and Implicit Schemes in the Simulation of the Heat
Equation, \textbf{Appl.Sci.} 2021, 11, 4468.
https://doi.org/10.3390/app11104468.

\bibitem{gohmajidismail} J. Goh, A. A. Majid and A. I. Ismail, Cubic
B-Spline Collocation Method for One-Dimensional Heat and Advection-Diffusion
Equations, \textbf{Journal of Applied Mathematics}, Volume 2012, Article ID
458701, 8 pages doi:10.1155/2012/458701

\bibitem{kaskar} N. F. Kaskar, Modified Implicit Method for Solving One
Dimensional Heat Equation, \textbf{International Journal of Engineering
Research in Computer Science and Engineering}, 8(9) 2021, 1-6.

\bibitem{cruzmercedesribeiro} G. Lozande-Cruz, C.E. Rubio-Mercedes and J.
Rodrigues-Riberio, Numerical Solution of Heat Equation with Singular Robin
Boundary Condition, \textbf{Tendencias em Matematica Aplicada e Computacional%
}, 19(2), (2018), 209-220

\bibitem{hooshmandasl} M.R. Hooshmandasl, M.H. Heydari and F.M. Maalek
Ghaini, Numerical Solution of the One-Dimensional Heat Equation by Using
Chebyshev Wavelets Method, \textbf{Journal of Applied Computational
Mathematics}, 1(6) (2012) 1-7. DOI: 10.4172/2168-9679.1000122

\bibitem{sonek1} H. Sun and J. Zhang, A high-order compact boundary value
method for solving one-dimensional heat equations, \textbf{Numerical Methods
for Partial Differential Equations}, 19 (6) (2003) 846-857.

\bibitem{sonek2} N. Patel, J. U. Pandya, One-Dimensional Heat Equation
Subject to Both Neumann and Dirichlet Initial Boundary Conditions and Used A
Spline Collocation Method, \textbf{Kalpa Publications in Computing}, 2
(2017) 107-112.

\bibitem{sonek3} T. Tarmizi, E. Safitri, S. Munzir and M. Ramli, On the
numerical solutions of a one-dimensional heat equation: Spectral and Crank
Nicolson method, \textbf{AIP Conference Proceedings} 2268, 2020.

\bibitem{8} C. de Boor, "A practicle guide to splines", \textbf{Applied
Mathematical Sciences}, (2001).

\bibitem{9} C. de Boor, "On calculating with B-splines", Journal of
Approximation Theory, (1972) 50-62.

\bibitem{14} A.R. Bahad\i r, Application of Cubic B-spline Finite Element
Technique to the Thermistor Problem, \textbf{Applied Mathematics and
Computation}, 149(2), (2004) 379-387.

\bibitem{15} B. Saka, I. Dag, Quartic B-spline Collocation Method to the
Numerical Solutions of the Burgers' Equation, \textbf{Chaos, Solitons
Fractals}, 32(3) (2007) 1125-1137.

\bibitem{16} B. Saka, I. Dag, A. Boz, B-spline Galerkin Methods for
Numerical Solutions of the Burgers' Equation, \textbf{Applied. Mathematics.
and Computation}., 166(3), (2005) 506-522.

\bibitem{17} M.A. Ramadan, T.S El-Danaf, F.E.I. Abd Alaal, A Numerical
Solution of the Burgers' Equation Using Septic B-splines, \textbf{Chaos,
Solitons Fractals}, 26(3) (2005) 795-804.

\bibitem{20} A.K. Khalifa, K.R. Raslana, H.M. Alzubaidi, A Collocation
Method with Cubic B-splines for Solving the MRLW Equation, J\textbf{ournal
of computational and Applied Mathematics} 212 (2008) 406-418.

\bibitem{21} A.H.A. Ali, L.R.T. Gardner, G.A. Gardner, A collocation
solution for Burgers' Equation Using Cubic Spline finite elements,\textbf{\
Computer Methods in Applied Mechanics and Engineering}., 100(3) (1992)
325-337.

\bibitem{hermiteek3} A. Yosaf, S.U. Rehman, F. Ahmad, M.Z. Ullah, A.S.
Alshomrani, Eight-Order Compact Finite Difference Scheme for 1D Heat
Conduction Equation, Advances in N\"{u}merical Analysis, (2016)

\bibitem{hermiteek4} F. Han, W. Dai, New higher-order compact finite
difference schemes for 1D heat conduction equations, Applied Mathematical
Modelling, 37 (2013) 7940-7952

\bibitem{newtrends} S. Kutluay, N. M. Ya\u{g}murlu and A. S. Karaka\c{s}, An
Effective Numerical Approach Based on Cubic Hermite B-spline Collocation
Method for Solving the 1D Heat Conduction Equation, \textbf{New Trends in
Mathematical Sciences, }10, No. 4, 20-31 (2022)
doi:http://dx.doi.org/10.20852/ntmsci.2022.485

\bibitem{septicek1} A. Kumari, V. K. Kukreja, Shishkin mesh based septic
Hermite interpolation algorithm for time-dependent singularly perturbed
convection--diffusion models, J\textbf{ournal of Mathematical Chemistry},
2022) 60:2029--2053 doi:https://doi.org/10.1007/s10910-022-01399-8
\end{thebibliography}
\end{document}